# Optimal control of multiple Markov switching stochastic system with application to portfolio decision

Jianmin Shi

*Haitong Securities Co. Ltd., Shanghai, China.*
*School of mathematics and statistics, Wuhan University, Hubei, China;*

**Abstract** In this paper we set up an optimal control framework for a hybrid stochastic system with dual or multiple Markov switching diffusion processes, while Markov chains governing these switching diffusions are not identical as assumed by the existing literature. As an application and illustration of this model, we solve a portfolio choice problem for an investor facing financial and labor markets that are both regime switching. In continuous time context we combine two separate Markov chains into one synthetic Markov chain and derive its corresponding generator matrix, then state the HJB equations for the optimal control problem with the newly synthesized Markov switching diffusion. Furthermore, we derive explicit solutions and value functions under some reasonable specifications.

Address correspondence to Jianmin Shi, 689 Guangdong Rd., Huangpu District, Haitong Securities Co., Ltd., Shanghai, CHINA; e-mail: shijm@htsec.com.

# 1.Introduction

In this paper, we will consider optimal control problems subject to an underlying stochastic system with dual or multiple regime switching diffusions following different Markov chains, as illustrated by an application to portfolio decision under regime switching financial and labor market.

In a real world, we learn that the properties of many natural or social systems is nonlinear or semi-linear rather than pure linear. For example, many practical systems may experience abrupt changes in their structure and parameters caused by such phenomena as component failures or repairs, changing subsystem interconnections in an electric control system, or abrupt environmental disturbances in an ecological or economic system, etc. An appropriate way to express this dynamics is stochastic hybrid control systems, which combines a part of the state that is driven by classic Brownian motion and another part of the state that takes discrete values. One of the important classes of the hybrid systems is the stochastic differential equations(SDE) with Markov switching. In this setting, the system will switch from one mode or state to another in a random way, and the switching between the modes is governed by a Markov chain. Optimal control or dynamic programming for such a hybrid stochastic system has been studied in different application fields(e.g. [39] [9] [45]).

A prominent feature of these existing research is that the structural or regime changes of the state variable dynamics are assumed to be single-sourcing, and therefore regime or structural changes of all parameters are uniformly represented by one common Markov chain. However, in real world, the sources of state or regime changes of stochastic systems are usually not the same but diverse. For instance, in an ecological system the size of prey populations is influenced by both availability of resources and predation pressure, which both show cyclical change pattern([40][40]). Obviously, characterizing these two cyclical changes of different sources by just one same Markov chain is unreasonable, though they may be interacted. Similarly, in an economic control system, the evolution of state variables (e.g. household's wealth, producer's profits etc.) depends on both business cycles and financial cycles, which have different generating mechanisms and are unsuitable to be represented by one common Markov chain (about interaction and differences between business and financial cycle, see [11] or [6]). Despite its significance from both theoretic and practical point of view, relaxing the assumption that all structural or regime changes of stochastic control system follow one same Markov chain is not a trivial or marginal task but technically challenging. Directly expanding the number of Markov chains in the control system will complicate the related optimization problem substantially and make existing dynamic programming techniques intractable.[1] In order to circumvent this difficulty, in this paper, we will seek a new way to set up a flexible and analytically tractable optimal control framework for stochastic systems with dual or multiple Markov-modulated diffusions, applying to an portfolio optimization problem with both regime switching financial and labor markets.

---

[1] Someone may think it's unnecessary or easy combining two or more separate Markov chains into a single chain, if a single Markov chain with an enlarged state space can play an equivalent role or two Markov chains can be combined into a single chain just by taking the Kronecker product of these processes. However, these methods still leave core problems unsolved, such as how to get the enlarged or combined chain's transition rates or generators out of its initial component chains.

Classic portfolio optimization model in continuous-time context was firstly established by the seminal work of Merton ([30][31]). In recent years one significant generalization of Merton model was characterizing structural or cyclical state changes driving by relative long-term factors in the financial markets with regime switching model originally proposed by Hamilton ([18][19]). Maximization of expected utility from consumption and/or terminal wealth in financial markets with regime switching has been studied, such as [24][47][39][27] etc. Another significant adaption of the Merton model to the real life is introducing stochastic stream of labor income as an important flow of the agent's wealth in additional to an initial endowment. [22] find that wage income is the most important source of wealth for the typical household. [13][26][23] and [32] successively propose the portfolio choice problem for an investor facing imperfectly hedgeable stochastic income. The above two significant generalizations have greatly enriched the classic Merton model. However, up to our knowledge, these generalizations are developed isolately. The regime switching models generally ignores the effect of stochastic labor income or labor markets. Similarly, the optimal investment models incorporating the labor income usually overlook the cyclical or structural regime changes of financial markets. In fact, not only financial markets but also labor markets are characterized by such a regime switching feature([10] [41] [36]). Therefore our aforementioned proposal of optimal control framework for dual or multiple Markov-modulated stochastic systems is well suited for solving an investor's portfolio optimization problem under two or more different regime switching markets.

In summary, in this paper we will set up a finite horizon optimal control framework for stochastic systems with dual or multiple different Markov-modulated diffusions and apply it to a portfolio optimization problem with dual regime-switching dynamics, i.e. how an agent who receiving income from a regime-switching labor market dynamically allocate her wealth under a regime-switching financial market. The rest of this paper is organized as follows. In the section 2, we present the optimization problem and formulate the model we are working with. Section 3 develops the methodology and presents the associated HJB equations. The application to portfolio choice is constructed in section 4, including explicit solution and related economic analysis of the results. Section 5 concludes.

**2. Problem formulation and model setup**

Through this paper we denote by $T > 0$ a deterministic finite horizon and let $(\Omega, \mathcal{F}, P)$ be a complete filtered probability space with filtration $\mathcal{F} = \{\mathcal{F}_t : t \in [0, T]\}$ satisfying the usual conditions, i.e., $\mathcal{F}$ is an increasing, right-continuous filtration and $\mathcal{F}_0$ contains all P-null sets. Now we consider the following stochastic processes: (i) a m-dimensional standard Brownian motion $\{B_t : 0 \leq t \leq T\}$ defined on the probability space $(\Omega, \mathcal{F}, P)$. (ii) two continuous time stationary Markov chains $\{\varepsilon(t): 0 \leq t \leq T\}$ and $\{\zeta(t): 0 \leq t \leq T\}$ with finite space states $\mathcal{S}^\varepsilon = \{1, 2, \dots, m\}$ and $\mathcal{S}^\zeta = \{1, 2, \dots, n\}$ respectively, where *m* and *n* are the number of regimes. Assume that the Markov chain ε(t) has strongly irreducible generator $Q^\varepsilon = [q_{ij}^\varepsilon]_{m \times m}$, where $q_{ii}^\varepsilon := -\lambda_i^\varepsilon < 0$ and $\sum_{j \in S} q_{ij}^\varepsilon = 0$ for every regime $i \in S^\varepsilon = \{1, 2, \ldots, m\}$; the Markov chain ζ(t) has a strongly irreducible generator $Q^\zeta = [q_{ij}^\zeta]_{n \times n}$, where $q_{ii}^\zeta := -\lambda_i^\zeta < 0$ and $\sum_{j \in S} q_{ij}^\zeta = 0$ for every regime $i \in S^\zeta = \{1, 2, \ldots, n\}$. We assume that the Markov processes ε and ζ are independent of the Brownian motion B.

For our aim, we firstly make the following standard technical assumptions:

**Assumption 2.1.** *The maps $f: [0,T] \times \mathcal{R}^N \times U \times S \times S \to \mathcal{R}^N$, $g: [0,T] \times \mathcal{R}^N \times U \times S \times S \to \mathcal{R}^{N \times M}$, $\Phi: [0,T] \times \mathcal{R}^N \times U \times S \times S \to \mathcal{R}^N$, $\Psi: [0,T] \times \mathcal{R}^N \times U \times S \times S \to \mathcal{R}^N$ are such that: (i) for each fixed $i, j \in S$, $f(.,.,.,i,j), g(.,.,.,i,j), \Phi(.,.,i,j), \Psi(.,.,i,j)$ are uniformly continuous; (ii) for each fixed $i, j \in S$, there exists $K > 0$ such that for $\varphi(t,x,u,i,j) = f(t,x,u,i,j), g(t,x,u,i,j), \Phi(t,x,u,i,j), \Psi(t,x,u,i,j)$, we have*
$$|\varphi(t,x,u) - \varphi(t,y,u)|^2 < K|x-y|^2, \quad |\varphi(t,0,u)|^2 < K.$$

Let U be a separate metric space, for a control variable $u$, we consider a stochastic controlled system with dual Markovian switching, for $t \in [0, T]$

$$dX(t) = f(t, X(t), u(t), \varepsilon(t), \zeta(t))dt + g(t, X(t), u(t), \varepsilon(t), \zeta(t))dB_t \quad (2.1)$$

$$X(0) = x, \varepsilon(0) = i, \zeta(0) = j$$

together with the following objective functional

$$J(t, x, i, j; u) = E[\int_0^T \Phi(t, X(t), u(t), \varepsilon(t), \zeta(t))dt + \Psi(T, X(T), \varepsilon(T), \zeta(T))] \quad (2.2)$$

where $(X(.), \varepsilon(.), \zeta(.)) \in \mathcal{R}^N \times S \times S$ denotes the state trajectory associated with a control trajectory $u(.)$ and starting from $(x, i, j)$ when $t=0$.

The control process $u: [0,T] \times \Omega \to U$ is said to be a admissible control if $u$ is measurable and $\{\mathcal{F}_t\}$-adapted, the stochastic differential equation (2.1) has a unique solution and

$$E[\int_0^T |\Phi(t, X(t), u(t), \varepsilon(t), \zeta(t))|dt] < \infty, \quad E[|\Psi(T, X(T), \varepsilon(T), \zeta(T))|] < \infty.$$

We denote the set f all admissible controls by $\mathcal{U}$.

Under the assumption 2.1, for any $u(.) \in \mathcal{U}$, Eq.(2.1) admits a unique solution(see, e.g.[35][15][2]), and the objective functional in Eq.(2.2) is well-defined. So the optimal control problem can be stated as follows:

**Problem 2.2.** *Select an admissible control $\hat{u} \in \mathcal{U}$ that maximizes $J(t, x, i, j; u)$ and find a value function V defined by $V(t, x, i, j) = \sup_{u \in \mathcal{U}} J(t, x, i, j; u)$ (2.3). The control $\hat{u}$ is regarded as an optimal control.*

The above optimal control problem is specified with particular initial time $t=0$. We can proceed to generalize it to any initial time $s \in [0, T]$. For any $(s, x, i, j) \in [0, T] \times \mathcal{R}^N \times S \times S$, consider the state equation:

$$dX(t) = f(t, X(t), u(t), \varepsilon(t), \zeta(t))dt + g(t, X(t), u(t), \varepsilon(t), \zeta(t))dB_t \quad (2.4)$$

$X(s) = x, \varepsilon(s) = i, \zeta(s) = j$

The corresponding objective functional is

$$J(s,x,i,j;u) = E[\int_s^T \Phi(t, X(t), u(t), \varepsilon(t), \zeta(t))dt + \Psi(T, X(T), \varepsilon(T), \zeta(T))] \quad (2.5)$$

where $(X(.), \varepsilon(.), \zeta(.)) \in \mathcal{R}^N \times S \times S$ denotes the state trajectory associated with a control trajectory $u(.)$ and starting from $(x, i, j)$ when $t$=s.

The value function is defined by $V(s,x,i,j) = \sup_{u \in \mathcal{U}} J(s,x,i,j;u)$ (2.6).

In the next section, we will discuss how to apply the dynamic programming principle for the above value function V and obtain the corresponding Hamilton-Jacobi-Bellman(HJB) equations associated with the optimal control problem under consideration.

### 3.Dynamic programming principle and HJB equation

In the existing literature(e.g.[39][2][27]), how to apply dynamic programming(DP) principle and derive the associated HJB equation to a stochastic system with state variable dynamics given by Markov switching SDE has been studied. All these studies assume that the structural or cyclical switching of the system or model parameters is governed by one single Markov chain. However, in our optimal control problem in above section(Eq. 2.1-3 or Eq. 2.4-6) switching SDE of state variable dynamics is governed by two different Markov chains. Because of complicated interactive relationship between different switching sources and curse of dimensionality, directly extending the optimal control framework of stochastic system governed by only one Markov chain to a dual or multiple Markov-modulated system is technically unpractical. In this followings, instead of applying mechanically the existing DP toolkits for single Markov-modulated systems, we will firstly combine two (or more) individual Markov chains into one compound Markov chain in continuous time context. Based on the newly combined continuous time Markov chain, we can then tackle with the relevant optimal control problem without obstacle.

Suppose that within a stochastic control system structural or regime changes of the state variable dynamics is not single-sourcing but stem from two different sources, which can be represented by two distinct finite state continuous time Markov chain $\varepsilon_t$ and $\zeta_t$ respectively. Without loss of generality, we assume both Markov chains take 2 states: $i, j \in \mathcal{S} = \{0,1\}$. Referring to the discrete time model of [20] ,we can construct a new 4-state continuous time Markov chain $\xi(t)$ from two 2-state continuous time Markov chains $\varepsilon(t)$ and $\zeta(t)$ as follows:

$\xi(t) = 1$, if $\varepsilon(t) = 0$ and $\zeta(t)$=0
$\xi(t) = 2$, if $\varepsilon(t) = 1$ and $\zeta(t)$=0
$\xi(t) = 3$, if $\varepsilon(t) = 0$ and $\zeta(t)$=1
$\xi(t) = 4$, if $\varepsilon(t) = 1$ and $\zeta(t)$=1

In the discrete time version, [20] show that the transition matrix of the new Markov chain composed of two independent Markov chains can be written as follows:

$$P = \begin{bmatrix} p_{11} & p_{21} & p_{31} & p_{41} \\ p_{12} & p_{22} & p_{32} & p_{42} \\ p_{13} & p_{23} & p_{33} & p_{43} \\ p_{14} & p_{24} & p_{34} & p_{44} \end{bmatrix} = \begin{bmatrix} p_{00}^{\varepsilon}p_{00}^{\zeta} & p_{10}^{\varepsilon}p_{00}^{\zeta} & p_{00}^{\varepsilon}p_{10}^{\zeta} & p_{10}^{\varepsilon}p_{10}^{\zeta} \\ p_{01}^{\varepsilon}p_{00}^{\zeta} & p_{11}^{\varepsilon}p_{00}^{\zeta} & p_{01}^{\varepsilon}p_{10}^{\zeta} & p_{11}^{\varepsilon}p_{10}^{\zeta} \\ p_{00}^{\varepsilon}p_{01}^{\zeta} & p_{10}^{\varepsilon}p_{01}^{\zeta} & p_{00}^{\varepsilon}p_{11}^{\zeta} & p_{10}^{\varepsilon}p_{11}^{\zeta} \\ p_{01}^{\varepsilon}p_{01}^{\zeta} & p_{11}^{\varepsilon}p_{01}^{\zeta} & p_{01}^{\varepsilon}p_{11}^{\zeta} & p_{11}^{\varepsilon}p_{11}^{\zeta} \end{bmatrix} \quad (3.1)$$

Where $p^{\varepsilon}$, $p^{\zeta}$, $p$ are the transition matrix of discrete Markov chain $\varepsilon(t)$, $\zeta(t)$, $\xi(t)$ respectively.

Under the continuous time circumstance, the corresponding core parameters, which should be derived from those of the original Markov chain $\varepsilon(t)$ and $\zeta(t)$, are the generator matrix Q of the newly synthesized or combined 4-state Markov chain $\xi(t)$:

$$Q = \begin{bmatrix} q_{11} & q_{21} & q_{31} & q_{41} \\ q_{12} & q_{22} & q_{32} & q_{42} \\ q_{13} & q_{23} & q_{33} & q_{43} \\ q_{14} & q_{24} & q_{34} & q_{44} \end{bmatrix} = \begin{bmatrix} -\lambda_1 & \lambda_1 p_{21} & \lambda_1 p_{31} & \lambda_1 p_{41} \\ \lambda_2 p_{12} & -\lambda_2 & \lambda_2 p_{32} & \lambda_2 p_{42} \\ \lambda_3 p_{13} & \lambda_3 p_{23} & -\lambda_3 & \lambda_3 p_{43} \\ \lambda_4 p_{14} & \lambda_4 p_{24} & \lambda_4 p_{34} & -\lambda_4 \end{bmatrix} \quad (3.2)$$

Here $\lambda_i$ is the transition rate of state or regime $i$ (i.e. the parameter of the exponential holding time for state or regime $i$) for Markov chain $\xi(t)$. $p_{ji}$ denotes the transition probability of embedded discrete-time Marko chain from regime $i$ to regime $j$ for $\xi(t)$.

In general, we may seek for the generator of the new compound Markov chain from its components via copula methods. A copula is a function that represents the joint distribution in terms of its marginals, and hence can be used to couple any discrete and/or continuous distributions. Following celebrated Sklar's theorem, Let *X, Y* be continuous r.v.'s, and let $F_{XY}$ be their joint distribution function with marginals $F_X$ and $F_Y$. Then there exists a unique copula C such that $F_{XY}(x,y) = C(F_X(x), F_Y(y))$ for all real *x, y*. Conversely, if C is a copula and $F_X$ and $F_Y$ are distribution functions, then $F_{XY}$ is a joint distribution function with marginals $F_X$ and $F_Y$. For detailed and standard reference of copulas, see [25]; [33].

Conventionally we construct a continuous time Markov chain(CTMC) with infinitesimal generator $Q = (q_{ij})$ as following: at each state $i$, wait for a sample time of an exponential distribution with rate $-q_{ii}$, and then jump to the next state according to the transition probability matrix *P* of its embedded discrete time Markov chain(DTMC), which can be computed out of *Q*: $P = -diag(Q)^{-1}Q + I$. Alternatively, a continuous time Markov chain(CTMC) with *m* states or regimes and generator matrix $Q = (q_{ij})$ can also be constructed as independent *m(m-1)* Poisson processes, each with parameters $\{q_{ij}, 1 \le i,j \le m\}$ corresponding to the events that the current state is *j* meanwhile the last state is *i* respectively, i.e. jumps from state *i* to state *j* (see [8];[21]). Based on the latter way of construction, we can combine two separate CTMCs by choosing a suitable copula that has Poisson distributions as its marginals. Gaussian copulas are a natural choice for connecting individual marginal distributions of counting data processes, such as binomial or Poisson distributions ([34]). A Gaussian copula distribution with Poisson marginal distributions is defined simply as:

$$G(x_1, x_2, \ldots, x_d|\theta) = C_\theta(F_1(x_1|\lambda_1), \ldots, F_d(x_d|\lambda_d))$$

Where $G$ represents a multivariate standard normal cumulative distribution function (CDF), $\theta$ denotes copula parameters, $F_i(x_i|\lambda_i)$ is the Poisson cumulative distribution function with mean parameter $\lambda_i$, here equal to corresponding transition rates of individual univariate CTMCs before combination. Differentiating $G(x_1, x_2, \ldots, x_d|\theta)$ gives the joint density or likelihood function with a correlation matrix $K$ as follows:

$$f(x) = \frac{1}{|K|^{\frac{1}{2}}} exp\left\{-\frac{1}{2}Z(K^{-1} - I)Z^T\right\} \prod_{i=1}^{d} \frac{1}{\sigma_i}\varphi(z_i)$$

where $Z = (z_1, z_2, \ldots, z_d)$, $z_i = \Phi^{-1}[F_i(x_i|\lambda_i)]$, $F_i(x_i|\lambda_i) \sim Poisson(\lambda_i)$. $\varphi$ denotes standard normal PDF, $\Phi$ represents a multivariate standard normal CDF. Here the density of copula is defined as: $(x) = \frac{1}{|K|^{\frac{1}{2}}} exp\left\{-\frac{1}{2}Z(K^{-1} - I)Z^T\right\}$.

When the marginal distributions $F_i(x_i|\lambda_i)$ are discrete such as Poisson distributions, we can obtain joint probability mass function using Gaussian copula as follows([38]):

$$g(y_1, \ldots, y_n; \Lambda, \Sigma) = P(Y_1 = y_1, \ldots, Y_n = y_n|\Lambda, \Sigma)$$

$$= \sum_{j_1=1}^{2} \sum_{j_2=1}^{2} \ldots \sum_{j_n=1}^{2} (-1)^{j_1+\cdots+j_n} C(u_{1j_1}, u_{1j_1}, \ldots, u_{1j_1}; \Sigma)$$

where $\Lambda = (\lambda_i, \lambda_i, \ldots \lambda_i)$, $C(u_{1j_1}, u_{1j_1}, \ldots, u_{1j_1}; \Sigma)$ is the Gaussian copula with association matrix $\Sigma$. $u_{i1} = F_i(y_i)$, $u_{i2} = F_i(y_i-)$, here $F_i(y_i-)$ is the left-hand limit of $F_i$ at $y_i$, for count-valued random variables, $F_i(y_i-) = F_i(y_i - 1)$. Taking derivative of joint probability mass function at time 0 we will obtain the transition rates or generator matrix of the newly combined Markov chain.

In case of independence between $\varepsilon(t)$ and $\zeta(t)$, taking advantage of the memoryless property of Poisson process, we can yield the elements of generator matrix of newly synthesized Markov chain $\xi(t)$ analytically without numerical estimation:

$q_{21}$
$$= \lim_{h \downarrow 0} \frac{P\left(\frac{\varepsilon(t+h) \neq 0, \zeta(t+h) = 0}{\varepsilon(t) = 0, \zeta(t) = 0}\right)}{h} = \lim_{h \downarrow 0} \frac{P((\varepsilon(t+h) \neq 0)/(\varepsilon(t) = 0)P((\zeta(t+h) = 0/(\zeta(t) = 0))}{h} = \lim_{h \downarrow 0} \frac{(\lambda_1^\varepsilon h + o(h))(1 - \lambda_1^\zeta h + o(h))}{h}$$
$$= \lambda_1^\varepsilon$$

$q_{31}$
$$= \lim_{h \downarrow 0} \frac{P\left(\frac{\varepsilon(t+h) = 0, \zeta(t+h) \neq 0}{\varepsilon(t) = 0, \zeta(t) = 0}\right)}{h} = \lim_{h \downarrow 0} \frac{P((\varepsilon(t+h) = 0)/(\varepsilon(t) = 0)P((\zeta(t+h) \neq 0/(\zeta(t) = 0))}{h} = \lim_{h \downarrow 0} \frac{(1 - \lambda_1^\varepsilon h + o(h))(\lambda_1^\zeta h + o(h))}{h}$$
$$= \lambda_1^\zeta$$

$$q_{41} = \lim_{h\downarrow 0} \frac{P\left(\genfrac{}{}{0pt}{}{\varepsilon(t+h)\neq 0,\zeta(t+h)\neq 0}{\varepsilon(t)=0,\zeta(t)=0}\right)}{h} = \lim_{h\downarrow 0} \frac{P((\varepsilon(t+h)\neq 0)/(\varepsilon(t)=0)P((\zeta(t+h)\neq 0/(\zeta(t)=0))}{h} = \lim_{h\downarrow 0} \frac{(\lambda_1^\varepsilon h + o(h))(\lambda_1^\zeta h + o(h))}{h}$$
$$= 0$$

$$q_{11} = -(q_{21} + q_{31} + q_{41}) = -(\lambda_1^\varepsilon + \lambda_1^\zeta) \quad (3.4)$$

Here $\lambda_i^\varepsilon$ and $\lambda_i^\zeta$ denote the transition rates of state or regime $i$ for Markov chain $\varepsilon(t)$ and $\zeta(t)$ respectively.

Following the same derivation, we can obtain $q_{21}, q_{22}, \ldots q_{44}$ similarly and therefore the whole generator matrix Q:

$$Q = \begin{bmatrix} q_{11} & q_{21} & q_{31} & q_{41} \\ q_{12} & q_{22} & q_{32} & q_{42} \\ q_{13} & q_{23} & q_{33} & q_{43} \\ q_{14} & q_{24} & q_{34} & q_{44} \end{bmatrix} = \begin{bmatrix} -(\lambda_1^\varepsilon + \lambda_1^\zeta) & \lambda_1^\varepsilon & \lambda_1^\zeta & 0 \\ \lambda_2^\varepsilon & -(\lambda_2^\varepsilon + \lambda_1^\zeta) & 0 & \lambda_1^\zeta \\ \lambda_2^\zeta & 0 & -(\lambda_1^\varepsilon + \lambda_2^\zeta) & \lambda_1^\varepsilon \\ 0 & \lambda_2^\zeta & \lambda_2^\varepsilon & -(\lambda_2^\varepsilon + \lambda_2^\zeta) \end{bmatrix} \quad (3.5)$$

Now we have a newly synthesized Markov chain $\xi(t)$, which takes 4 states $i \in \mathcal{S} = \{1,2,3,4\}$ and has strongly irreducible generator $Q = [q_{ij}]_{4\times 4}$, where $q_{ii} := -\lambda_i < 0$ and $\sum_{j\in\mathcal{S}} q_{ij} = 0$ for every regime $i \in \mathcal{S} = \{1,2,3,4\}$. As demonstrated above, the generator of the compound Markov chain $\xi(t)$ is derived from those of two original Markov chains composing of $\xi(t)$. With the newly synthesized Markov chain $\xi(t)$, we can reshape the optimal control problem(2.1-3) or (2.4-6) as follows:

For any $(s, x, i) \in [0, T] \times \mathcal{R}^N \times S$, consider the state equation:

$$dX(t) = f(t, X(t), u(t), \xi(t))dt + g(t, X(t), u(t), \xi(t))dB_t \quad (3.6)$$

$X(s) = x, \xi(s) = i$

Where $f = (f_1, \ldots, f_N)^T : [0, T] \times \mathcal{R}^N \times U \times S \to \mathcal{R}^N$, $g = (g_{kl})_{N\times M} : [0, T] \times \mathcal{R}^N \times U \times S \to \mathcal{R}^{N\times M}$.

The corresponding objective functional is

$$J(s, x, i; u) = E\left[\int_s^T \Phi(t, X(t), u(t), \xi(t))dt + \Psi(T, X(T), \xi(T))\right] \quad (3.7)$$

where $(X(.), \xi(.)) \in \mathcal{R}^N \times S$ denotes the state trajectory associated with a control trajectory $u(.)$ and starting from $(x, i)$ when $t=s$.

The optimal control problem under consideration is then stated as: For any $(s, x, i) \in [0, T] \times \mathcal{R}^N \times S$, find $\bar{u} \in \mathcal{U}$ such that
$$J(s, x, i; \bar{u}) = \sup_{u\in\mathcal{U}} J(s, x, i; u) \quad (3.8)$$
and the corresponding value function is defined by
$$V(s, x, i) = \sup_{u\in\mathcal{U}} J(s, x, i; u) \quad (3.9).$$

In the following, we will deduce the Hamilton-Jacobi-Bellman(HJB) equation associated with the above dynamic optimization problem. The HJB approach will apply the so-called Bellman optimality principle, which states that an optimal policy has the property that whatever the initial state and initial decision are, the remaining decisions must constitute an optimal policy with regard to the state resulting from the first decision ([4]). This is formally expressed in the following equation:

$$V(t,x,i) = \sup_{u \in \mathcal{U}} E[(V(\hat{t}, X(\hat{t}), \xi(\hat{t}))]; \hat{t} \in [t,T] \quad (3.10)$$

We will now use the principle of optimality to derive the corresponding HJB equation, a sequence of partial differential equation(PDE), whose solution is the value function of the optimal control problem under consideration here. Let $C^{1,2}([0,T] \times \mathcal{R}^N; \mathcal{R})$ denote all functions $V(t,x,i)$ on $[0,T] \times \mathcal{R}^N$ that are continuously differentiable in $t$, continuously twice differentiable in $x$ for each $i \in \mathcal{S}$. Corresponding to $V(.,i) \in C^{1,2}([0,T] \times \mathcal{R}^N; \mathcal{R})$ we firstly define an operator $L_i(u)V$ for each $i \in \mathcal{S}$ by

$$L_i(u)V(t,x,i) := \sum_{k=1}^{N} \frac{\partial V}{\partial x_k}(t,x,i) f_k(t,x,u,i)$$

$$+ \frac{1}{2} \sum_{k=1}^{N} \sum_{l=1}^{N} \frac{\partial^2 V}{\partial x_k x_l}(t,x,i) \sum_{q=1}^{M} g_{kq}(t,x,u,i) g_{lq}(t,x,u,i)$$

$$+ \frac{\partial V}{\partial t}(t,x,i) \quad (3.10)$$

In what follows, we give the equation of optimality for the optimal policy.

**Theorem 3.1** *For each $\in \mathcal{S}$, V(., i) is a solution of the following terminal problem of the Hamilton-Jacobi-Bellman(HJB) equation*

$$sup_{u \in \mathcal{U}}[L_i(u)V(t,x,i) + \Phi(t,x,u,i)] = \lambda_i V(t,x,i) - \sum_{j \in \mathcal{S} \setminus \{i\}} q_{ij} V(t,x,j)$$

(3.11)
*with the terminal condition $V(T,x,i) = \Psi(T,x,i)$.*
*Proof.* See Appendix.

## 4. Application to portfolio choice problem
### 4.1 regime switching financial and labor market set-up

In the application part, we will consider dynamic portfolio selection problem of a representative economic agent(an investor) who receives a stream of stochastic labor income and faces financial markets consisting of a riskless asset(bond) and a risk asset(stock). The information structure of the financial markets can be described by a probability space $(\Omega, \mathcal{F}, P)$, on which a standard Brown motion $B_t$ and a continuous-time finite-state Markov chain $\varepsilon_t$ are defined and we assume that they are independent. The filtration $\{\mathcal{F}_t\}_{t \geq 0}$ is generated by both the Brownian motion $B_t$ and the Markov chain $\varepsilon_t$.

In the financial markets, it is assumed that two assets are traded: a risky asset (which we refer to as stock) with price $P$, and a riskless bond with price $R_t = R_0 e^{rt}$ at time t where $r$ is a constant rate of interest. The stock price follows

$$dP(t) = P(t)(\alpha_{\epsilon(t)}dt + \sigma_{\epsilon(t)}dB_t) \quad (4.1)$$

with initial price $P(0) > 0$ and initial state $\varepsilon_0$.

Here $\alpha, \sigma$ is the expected rate of return and volatility of the stock respectively. B is a standard Brownian motion. $\varepsilon_t \in \varepsilon$, represents the state or regime of the stock market at time t. $\varepsilon$ is an continuous-time, stationary, finite-state Markov chain(for simplicity, we assume 2 states or regimes, 'bullish' and 'bearish', in the stock market, see [16]).Furthermore, we assume that the Markov chain has a strongly irreducible generator $Q^\epsilon = [q_{ij}^\epsilon]_{2 \times 2}$ ,where $q_{ii}^\epsilon := -\lambda_i^\epsilon < 0$ and $\sum_{j \in S} q_{ij}^\epsilon = 0$ for every regime $i \in S = \{0,1\}$.

In additional to invest in the financial markets, the agent also receives labor income over time, the income rate $Y_t$ following

$$dY(t) = \mu(Y_t, \zeta_t, t)dt + \delta(Y_t, \zeta_t, t)dZ_t \quad (4.2)$$

with initial price $Y(0) > 0$ and initial state $\zeta_0$. Similarly, here $\mu, \delta$ is the expected mean and volatility of the stochastic income rate respectively.

The correlation between $dB_t$ and $dZ_t$ is $\rho dt$, where $-1 \leq \rho \leq 1$. it is convenient to write Z as $Z_t = \rho B_t + \sqrt{1-\rho^2} W_t$ where W is a standard Brownian motion independent of B. $\zeta(t) \in \zeta$, represents the state or regime of the labor market at time t. $\zeta$ is an continuous-time, stationary, finite-state Markov chain(for simplicity, we assume 2 states or regimes, 'boom' and 'recession', in the labor market).Furthermore, we assume that the Markov chain has a strongly irreducible generator $Q^\zeta = [q_{ij}^\zeta]_{2 \times 2}$ ,where $q_{ii}^\zeta := -\lambda_i^\zeta < 0$ and $\sum_{j \in S} q_{ij}^\zeta = 0$ for every regime $i \in S = \{0,1\}$.

As discussed in section 3, we know that characterizing the regime switching financial markets and labor market by two different Markov chains respectively will complicate the related optimization problem substantially and make the dynamic programming technique in the existing literature intractable. So in the following we will firstly use the method described in section 3 to combine two different Markov chains into one compound Markov chain. With the newly compound continuous time Markov chain, we can then tackle with the relevant optimal control problem without a hitch.

Suppose two finite state continuous time Markov chain ,which represent the phases or regimes of financial market($\varepsilon_t$) and labor market($\zeta_t$)respectively. Without loss of generality, we assume both Markov chains take 2 states $\{0,1\}$, i.e. there are "bull" and "bear" two regimes in stock market, and the labor market shifts between 2 states of "boom" and "recession". Based on these setting, we can construct a new 4-state Markov chain $\xi(t)$ from $\varepsilon(t)$ and $\zeta(t)$ as follows:

$\xi(t) = 1, \text{if } \varepsilon(t) = 0 \text{ and } \zeta(t)=0$
$\xi(t) = 2, \text{if } \varepsilon(t) = 1 \text{ and } \zeta(t)=0$

$\xi(t) = 3$, if $\varepsilon(t) = 0$ and $\zeta(t)=1$

$\xi(t) = 4$, if $\varepsilon(t) = 1$ and $\zeta(t)=1$

The generator Q of the newly synthesized or combined 4-state Markov chain $\xi(t)$ can be derived by the approach as described in section 3.

Underlying the newly synthesized Markov chain $\xi(t)$, the stock price and labor income variable process in (4.1) and (4.2) will satisfy the following Markov-modulated stochastic differential equations:

$$dP(t) = \alpha_{\xi(t)}P(t)dt + \sigma_{\xi(t)}P(t)dB_t \quad (4.3)$$

$$dY(t) = \mu(Y_t, \xi_t, t)dt + \delta(Y_t, \xi_t, t)dZ_t \quad (4.4)$$

The agent chooses a portfolio $\pi = \{\pi(t), t \geq 0\}$, representing the fraction(not ratio) of wealth X invested in the risky asset. The fraction of wealth invested in the riskless asset at time $t \geq 0$ is then $X(t) - \pi(t)$. Here a technical condition need to be satisfied: A portfolio vector process is an uncertain stochastic vector process $\pi$ such that $E\left[\int_0^t \pi(s)\sigma_{\xi(t)}ds\right] < +\infty$ for all $t \geq 0$.

The agent's wealth process $X = \{X(t), t \geq 0\}$, determined by his specific portfolio holding and the inflow of stochastic income, can be written as

$$dX(t) = [rX(t) + \pi(t)(\alpha_{\xi(t)} - r) + Y(t)]dt + \pi(t)\sigma_{\xi(t)}dB_t \quad (4.6)$$

with initial wealth $X(0) > 0$ and initial state $\xi(0) \in S = \{1,2,3,4\}$.

Here it should be noted that we assume no consumption in the agent's wealth process. Its exclusion from optimal portfolio choice problem is not without precedent( see [46][23] etc.). Such consumption would not change the nature of our results. If consumption is involved, it is unlikely that we could obtain closed form solution to our problem. Portfolio choice models incorporating both consumption and incomplete income process have to be solved numerically( see [44]).

Next we introduce the utility functions $U(.)$ of terminal wealth, which are assumed to be twice differentiable, strictly increasing, and concave. Moreover, $U'(0) = \infty$, $U'(\infty) = 0$. It is easily to see that there exists a constant K such that $U(z) \leq K(1 + z^2)$.

For $t \in [0, T]$, we consider the problem of an investor with expected utility over terminal wealth which can be maximized by the selection of portfolio $\pi$. The objective function in defined as

$$J(t, x, y, i; \pi) := E[U(X(T))] \quad (4.7)$$

where for $s \in [t, T]$, the wealth X(s) process follows (4.6) with initial X(t)=x, Y(t)=y and initial state $\xi(t) = i$. The control process $\pi$ which is a portfolio process satisfying

$$E[U^-(X(T))] < +\infty \quad (4.8)$$

will be called an admissible control process. Here we denote $a^- := \max(0, -a)$ and $a^+ := \max(0, a)$ for every $a \in [-\infty, \infty]$. The set of all admissible controls will be denoted by $\bar{A}$.

Based on the objective function, the value function is defined by

$$V(t,x,y,i) = \sup_{\pi \in Ä} J(t,x,y,i;\pi) \qquad (4.9)$$

In summary, the agent wants to solve the following problem:

**Problem 4.1.** *Select an admissible control $\hat{\pi} \in Ä$ that maximizes (4.7) and find the value function defined by (4.9). The control $\hat{\pi}$ is called an optimal control or an optimal policy.*

## 4.2 The HJB Equation

In this subsection, we will derive the corresponding HJB equation, a sequence of partial differential equation(PDE), whose solution is the value function of the optimal control problem under consideration here. Let $C^{1,2}([0,T] \times \mathbb{R}^N; \mathbb{R})$ denote all functions $V(t,x,y,i)$ on $[0,T] \times \mathbb{R}^N$ that are continuously differentiable in t, continuously twice differentiable in x for each $i \in \mathcal{S}$. Corresponding to (3.10), we firstly define an operator $L_i(\pi)V$ for each $i \in \mathcal{S}$ by

$$L_i(\pi)V(t,x,y,i) := \frac{1}{2}\pi^2 \sigma_i^2 V_{xx} + (rx+y)V_x + \pi(\alpha_i - r)V_x + \frac{1}{2}\delta(y,i,t)^2 V_{yy} +$$

$$\mu(y,i,t)V_y + \pi \sigma_i \rho \delta(y,i,t) V_{xy} + V_t \quad (4.10)$$

Here subscripts on V denote partial derivatives, and we have suppressed the arguments of the functions for notational simplicity.

Applying theorem 3.1, we have the equation of optimality for the optimal policy as follows:

**Theorem 4.2** *For each $\in \mathcal{S}$, V(., i) is a solution of the following terminal problem of the Hamilton-Jacobi-Bellman(HJB) equation*

$$\sup_{\pi \in Ä}[L_i(\pi)V(t,x,y,i)] = \lambda_i V(t,x,y,i) - \sum_{j \in \mathcal{S} \setminus \{i\}} q_{ij} V(t,x,y,j) \quad (4.11)$$

*with the terminal condition $V(T,x,y,i) = U(x)$.*

Note that from (4.11), we will need to solve a system of coupled partial differential equations(PDE), instead of dealing with only one PDE as in the case with classical non-Markov-modulated diffusions.

Differentiating the left of (4.11) with respect to $\pi$ gives the first-order condition(FOCs):

$$\hat{\pi}(t) = -\frac{(\alpha_i - r)V_x + \sigma_i \delta(y,i,t) \rho V_{xy}}{\sigma_i^2 V_{xx}} = \frac{-(\alpha_i - r)V_x}{\sigma_i^2 V_{xx}} - \frac{\delta(y,i,t)\rho V_{xy}}{\sigma_i V_{xx}} \quad (4.12)$$

We can observe that the optimal portfolio in the risk asset(stock), $\hat{\pi}$ given in (4.12) is comprised of two components: the first part $\frac{-(\alpha_i - r)V_x}{\sigma_i^2 V_{xx}}$ is similar to the form of the classic Merton mean-variance portfolio strategy followed in the absence of stochastic income. The main difference is that stock returns $(\alpha_i)$ and volatilities$(\sigma_i)$ are shifting with different regimes of financial market, whereas in the Merton's classic model they are constant. Furthermore, the value function V here, which is also affected indirectly by the switching regimes, may not necessarily be identical to that of the Merton mean-variance strategy. The second part of our optimal portfolio strategy $\frac{-\delta(y,i,t)\rho V_{xy}}{\sigma_i V_{xx}}$ , in which $V_{xy}$ measures the sensitivity of the marginal utility of wealth to the stochastic income or the attitude to changes in stochastic income, can be interpreted as an intertemporal hedging against the variations of stochastic income. It depends on the volatilities of both financial market and the labor market as well as the value function V under different regimes. It's noteworthy that the hedging term will disappear if $\rho = 0$, in this case where the stock returns and income are uncorrelated the stochastic income can't be hedged with the traded asset. However, even in this case, the stochastic income still affect the optimal portfolio strategy via the value function.

The expression of optimal control $\hat{\pi}$ in (4.12) may partially shed some light on the intuitive interpretation of the optimal portfolio strategy. For example, one can speculate that the volatility of the risky asset and market price of risk under different market regimes should have impact on the optimal portfolio weights in stock, as these factors influence the attractiveness of the risk asset to the investor. Nevertheless, the representation in (4.12) indicates that thorough understandings of the optimal strategy require detailed knowledge about the value function, which will be the focus in the following.

*4.3 Construction of the solution and optimal strategies*

Given the HJB equation (4.11),in this section we will attempt to derive the solution and optimal policies for the stochastic control problem with regime-switching financial and labor market. As we know, for optimal control problem, analytical solutions allowing for robust comparative static analysis take advantage over the numerical solutions in many aspects. However, generally explicit solutions to dynamic optimization problems with nonlinear parameters are rarely available. As far as we know, a few papers([48][39][28]) study the portfolio selection problem in regime switching models and derived explicit solutions for the corresponding HJB equations and optimal policies. Even in these papers, the portfolio optimization models don't consider stochastic income, which is another nontrivial and uncertain driver of the investor's wealth besides financial assets. To the best our knowledge, our paper is the first attempt to solve the Markov-modulated optimal portfolio problem facing an investor with both stochastic financial revenue and income flow.

Let's start our journey by specifying the investor's utility function. The optimal strategy $\hat{\pi}$ characterized by (4.12) is expressed in terms of derivatives of the value function solving HJB equation, which generally depends on investor preferences and the investment horizon. As a result, solution of the optimal investment problem (4.9) is determined by these variables and the market parameters in an indirect way. So, to gain more insightful and interpretable results we need to specialize to a particular utility for the agent.

In this paper we will assume the negative exponential utility function

$$U(x) = -\left(\frac{1}{\gamma}\right)e^{-\gamma x}, \gamma > 0. \quad (4.13)$$

We prefer to this utility function because there are evidence from option prices that exponential utility provides a better representation of preferences than other forms of utility([6]). Given the negative exponential utility, we conjecture the following form for the value function

$$V(t, X_t, Y_t, i) = -\frac{1}{\gamma} e^{-\gamma X_t e^{r(T-t)}} g(t, Y_t, i), \text{with } g(T, y, i) = 1 \quad (4.14)$$

for some unknown function $g(t, y, i)$, which to be derived in the next step. Plugging the corresponding computed $V_x, V_y, V_{xx}, V_{xy}, V_{yy}$ together into the equation (4.12), we then get

$$\hat{\pi}(t) = \frac{(\alpha_i - r)}{\gamma \sigma_i^2 e^{r(T-t)}} + \frac{\delta(y,i,t)\rho g_y}{\gamma \sigma_i e^{r(T-t)} g} \quad (4.15)$$

As described earlier, the first part of $\hat{\pi}$ has the same structure as the Merton([30]) optimal portfolio under each regime $i$ as if stochastic income inflows were ignored. The second hedging term is independent of the investor's initial wealth but depend on the investor's risk aversion $\gamma$, the correlation between financial market and labor market $\rho$ and stochastic income. However, deeper understanding of the hedging term in (4.15) requires additional knowledge of the function $g$ or the value function V.

In order to get the solution function, we substitute equation(4.14) and (4.15) into equation (4.11) and obtain the following system of nonlinear partial differential equations(PDEs) of the value function $g$( or V )satisfying the HJB equation:

$$g_t(t, y, i) + g_y \left( \mu(y, i, t) - \frac{(\alpha_i - r)\delta(y,i,t)\rho}{\sigma_i} \right) + \frac{1}{2} g_{yy} \delta(y, i, t)^2 - g e^{r(T-t)} \gamma y -$$

$$\frac{(\alpha_i - r)^2}{2\sigma_i^2} g - \frac{g_y^2}{2g} \delta(y, i, t)^2 \rho^2 = -\sum_j q_{ij} g(t, y, j) \quad (4.16)$$

with the terminal condition $g(T, y, i) = 1$ (4.17)

Obviously the above systems (4.16-17) are highly non-linear and difficulty to obtain close form analytic solutions because of regime switching and the coupling of the equations. In fact for such a nonlinear and Markov modulated diffusion system even if finding numerical solutions is challenging. Up to our knowledge, in the literature there are few efficient numerical methods for stochastic control problem of regime switching systems. [37] may be an exception, which develops numerical algorithms for regime switching controlled diffusions and regime-switching jump diffusions. Even feasible, this numerical methods is still formidable because of computation intensive. Fortunately, under reasonable specifications some special forms of HJB equations in (4.16) can obtain explicit expressions that can be easily computed numerically. In the following we will present such special models. The first one is the case of $\rho = 0$, which means the Brownian motions of stock returns and stochastic income process are uncorrelated. Another special case doesn't put the constraint $\rho = 0$, but instead of assuming that the income are normally distributed under different regimes, which means the Markov-switching mean and volatility of the labor come process are independent of the income level.

### 4.3.1 Dual Markov-modulated models with independent Brownian motions

In this subsection, we set $\rho = 0$ in the general model(4.3-6),which states that the Brownian motions of stock returns and stochastic income process do not exhibit instantaneous correlation. However, this doesn't mean the two processes are independent at all, as they can be additionally correlated by the joint Markov chain $\xi(t)$.

It's evident that if $\rho = 0$, the system of PDEs (4.16) for function $g$ collapse to:

$$g_t(t,y,i) + g_y\mu(y,i,t) + \frac{1}{2}g_{yy}\delta(y,i,t)^2 - g\left[e^{r(T-t)}\gamma y + \frac{(\alpha_i-r)^2}{2\sigma_i^2}\right] = -\sum_j q_{ij}g(t,y,j)$$

or

$$g_t(t,y,i) + g_y\mu(y,i,t) + \frac{1}{2}g_{yy}\delta(y,i,t)^2 + \sum_j q_{ij}g(t,y,j)$$

$$- g\left[e^{r(T-t)}\gamma y + \frac{(\alpha_i - r)^2}{2\sigma_i^2}\right] = 0$$

(4.18)

with the terminal condition $g(T,y,i) = 1$ (4.19)

Under regular conditions for process Y and function $g$, applying the Feynman-Kac formula for Markov switching diffusions, the solution to (4.18) is given by:

$$g(t,y,i) = E^{y,i}\left[\exp\left(-\int_t^T (e^{r(T-s)}\gamma Y_s + \frac{(\alpha_{\xi(s)}-r)^2}{2\sigma_{\xi(s)}^2})ds\right)\right] \quad (4.20)$$

From (4.14), (4.15) and (4.20), we can obtain the following result:

**Proposition 4.3.** *The investor's value function for the optimal portfolio allocation problem 4.1 with dual Markov switching diffusion and independent Brownian motions is given by*

$$V(t,x,y,i) = -\frac{1}{\gamma}e^{-\gamma x e^{r(T-t)}} E^{y,i}\left[\exp\left(-\int_t^T (e^{r(T-s)}\gamma Y_s + \frac{(\alpha_{\xi(s)}-r)^2}{2\sigma_{\xi(s)}^2})ds\right)\right] \quad (4.21)$$

*and the corresponding optimal investment strategy is*

$$\hat{\pi}(t) = \frac{(\alpha_i-r)}{\gamma\sigma_i^2 e^{r(T-t)}} \quad (4.22)$$

*Proof.* See Appendix.

So in the case of $\rho = 0$ the optimal strategy consists only on the mean-variance portfolio just as the Merton investment strategy, the hedging term disappears and the income diffusion has no effect on the optimal portfolio. This mean-variance term, which is proportional to the risk premium and reciprocal of the coefficient of absolute risk aversion, is intuitively clear, as higher excess return of the risky asset lead to higher investment, whereas higher variance or risk reduces the position in the corresponding asset. Consistent with the Merton model, the optimal portfolio expressed in units of time $T$ money is constant in every regime of the market, which also does not depend on the investor's wealth. Although the optimal portfolio is constant under constant market coefficients with just one regime of the markets, the situation will change with the switching of market regimes and the optimal portfolio does depend on the regime of the market: it is positively proportional to the expected excess return on the stock in such regime $(\alpha_i - r)$, and inversely proportional to the variance of the stock return in such regime $\sigma_i^2$, which means that the investor prefers to invest a smaller fraction of wealth in the risky asset when the market presents weak conditions or "bearish" regime. If market conditions are strong or in "bullish" regime(means positive and large market risk premium $(\alpha_i - r)/\sigma_i^2$, the investor will assign a greater fraction of his wealth to invest in the stock. So once the markets change its regime, the optimal portfolio will vary and isn't constant any more. From equation (4.22), we see that when the investor has low risk tolerance(high risk aversion $\gamma$), the fraction of wealth invested in the stock is even smaller in each regime. However, no matter how the investor's risk attitudes change, the relation that $\hat{\pi}$ of "bearish" regime always less than that of "bullish" regime remain the same. In other words, every investor (independent of the level of risk aversion) should allocate a higher proportion of his wealth in a risky asset during a bull market than during a bear market.

Equation (4.22) indicates that the stochastic income cannot be hedged against in the case where the diffusions of stock returns and income are uncorrelated($\rho = 0$) and the hedging term disappears and the nontraded income has no effect on the optimal portfolio. Even so, however, from (4.21) we can see that the existence of the income still impacts on the value function.

*4.3.2 Normally distributed income process under different regimes*

In this special case, some constraint is put on the labor income process instead of the correlation between the Brown motions of the stock return process and income flow process. In the literature opinions about the distribution of income earnings process is inconclusive. Besides normal distribution, another popular point of view is log-normality like that of asset prices. Nevertheless, this assumption isn't indisputable. For example, as documented by [1] and [17], the earnings shocks experienced by US workers display important deviations from the assumptions of log-normality and independence from age and earnings realizations ([12]). Our assumption of normal distribution isn't only for tractability. In fact, prior research by [14], [42] and [43] also treats normally distributed income as return of non-traded asset, differing because their income is received at the investment horizon ([23])

Under our normal distribution assumption, the model parameters of income earnings process are specified as:
$dY(t) = \mu_i dt + \delta_i dZ_t$ (4.23)

where $\mu_i$ and $\delta_i$ are different constants under different regime respectively.

Further, in order to find the solution, we assume separability of function $g(t,y,i)$ in $y$ and $i$, i.e. the following ansatz is considered:

$g(t,y,i) = h(t,i)e^{M(t)y}$ (4.24)

for some function $h:[0,T] \times i \to \mathcal{R}$ and $M:[0,T] \to \mathcal{R}$.

From (4.24) we can derive the following expressions(subscripts denote corresponding partial derivatives):

$g_t(t,y,i) = h_t e^{M(t)y} + h e^{M(t)y} M_t y$ (4.25)

$g_y(t,y,i) = h e^{M(t)y} M$ (4.26)

$g_{yy}(t,y,i) = h e^{M(t)y} M^2$ (4.27)

Substituting above expressions (4.24),(4.25-27) into HJB equations (4.16), we can get the following system of equations:

$h_t + h\left\{M_t y + M\left[\mu_i - \frac{\rho(\alpha_i - r)}{\sigma_i}\delta_i\right] - e^{r(T-t)}\gamma y - \frac{1}{2}(\rho^2 - 1)\delta_i^2 M^2 - \frac{(\alpha_i - r)^2}{2\sigma_i^2}\right\} =$
$-\sum_j q_{ij} h(t,j)$ (4.28)

In order for (4.28) to hold for any $y \geq 0$, we assume that the undetermined $M(t)$ make the coefficient on $y$ to be zero, which imply the following ODE for $M(t)$:

$M_t - e^{r(T-t)}\gamma = 0, M(T) = 0$ (4.29)

Further (4.28) and (4.29) jointly lead to the following system of PDEs for function $h(t,i)$:

$h_t + h\left\{\frac{(\alpha_i - r)^2}{2\sigma_i^2} + M\left[\mu_i - \frac{\rho(\alpha_i - r)}{\sigma_i}\delta_i\right] - \frac{1}{2}(\rho^2 - 1)\delta_i^2 M^2\right\} =$
$-\sum_j q_{ij} h(t,j), h(T,i) = 1$ (4.30)

In general, equations (4.30) can't be solved in a closed form, however, we can provide the probabilistic representation for function $h(t,i)$. And then compute it in general numerically. Based on the application of the Feynman-Kac theorem with Markov switching, we can obtain the probabilistic representation for function $h(t,i)$ in (4.30) as follows:

$h(t,i) = E^i\left[\exp\left(\int_t^T \{\frac{(\alpha_i - r)^2}{2\sigma_i^2} + M\left[\mu_i - \frac{\rho(\alpha_i - r)}{\sigma_i}\delta_i\right] - \frac{1}{2}(\rho^2 - 1)\delta_i^2 M^2\}ds\right)\right]$ (4.31)

where

$M = -\frac{1}{r}e^{r(T-t)}\gamma + \frac{1}{r}\gamma = -\frac{\gamma}{r}(e^{r(T-t)} - 1)$ (4.32)

From (4.14),(4.15)and (4.24),(4.31),(4.32), we can obtain the following investor's value function and corresponding optimal investment strategy:

$V(t,X_t,Y_t,i) = -\frac{1}{\gamma}e^{-\gamma X_t e^{r(T-t)}}g(t,Y_t,i) = -\frac{1}{\gamma}e^{-\gamma X_t e^{r(T-t)}}h(t,i)e^{M(t)y}$

(4.33)

Where functions $h, M$, are given as in the equation (4.31) and (4.32) respectively.

$$\hat{\pi}(t) = \frac{(\alpha_i - r)}{\gamma \sigma_i^2 e^{r(T-t)}} + \frac{\delta(y_t, i, t)\rho g_y}{\gamma \sigma_i e^{r(T-t)} g} = \frac{(\alpha_i - r)}{\gamma \sigma_i^2 e^{r(T-t)}} + \frac{\delta_i \rho M}{\gamma \sigma_i e^{r(T-t)}}$$

$$= \frac{(\alpha_i - r)}{\gamma \sigma_i^2 e^{r(T-t)}} - \frac{\delta_i \rho [e^{r(T-t)} - 1]}{r \sigma_i e^{r(T-t)}}$$

(4.34)

In the following proposition we summarize what we have just derived:

**Proposition 4.4.** *The investor's value function for the optimal portfolio allocation problem 4.1 with dual Markov switching diffusion and normally distributed income process is given by*

$$V(t, X_t, Y_t, i) = -\frac{1}{\gamma} e^{-\gamma X_t e^{r(T-t)}} g(t, Y_t, i) = -\frac{1}{\gamma} e^{-\gamma X_t e^{r(T-t)}} h(t, i) e^{M(t)y} \quad (4.35)$$

*Where functions $h, M$, are given as*

$$M = -\frac{\gamma}{r}(e^{r(T-t)} - 1) \quad (4.36)$$

$$h(t, i) = E^i \left[ \exp\left( \int_t^T \{ \frac{(\alpha_i - r)^2}{2\sigma_i^2} + M \left[ \mu_i - \frac{\rho(\alpha_i - r)}{\sigma_i} \delta_i \right] - \frac{1}{2}(\rho^2 - 1)\delta_i^2 M^2 \} ds \right) \right] (4.37)$$

*and the corresponding optimal investment strategy is*

$$\hat{\pi}(t) = \frac{(\alpha_i - r)}{\gamma \sigma_i^2 e^{r(T-t)}} - \frac{\delta_i \rho [e^{r(T-t)} - 1]}{r \sigma_i e^{r(T-t)}} \quad (4.38)$$

(The proof is similar to that of **Proposition 4.3**)

Equation (4.34) or (4.38) gives optimal investment in the risky assets in this case of normally distributed income process. As mentioned before, the first term in the equation just corresponds to the classic Merton strategy or the standard mean-variance optimal portfolio without stochastic income. The excess return of the risky asset influences positively the positions in the corresponding risky asset. The contrary holds for the volatility term. The only point different from classic Merton model is the role the Markov chain plays in the portfolio selection, i.e. once the markets shift its regime, the optimal mean-variance portfolio will vary and isn't constant any more ,since the risky asset's excess return and volatility vary with switching regimes.

The second term of the optimal strategy is a portfolio that hedges against unexpected changes of stochastic labor income. Besides the riskless rate $r$ and the remaining time to the investment horizon $(T-t)$, the hedge demand depends on the ratio of volatility of labor market and financial market under the corresponding switching regime $\delta_i/\sigma_i$, the correlation between Brown motions of financial market and labor market $\rho$. The income hedge is not investor specific, and is not influenced by the investor's risk preference. This contrasts to the first mean-variance term which depends on the stock's risk premium and investor's risk aversion. The sign of the hedge demand is determined by the correlation between financial market and labor market $\rho$. The reason is intuitive. The hedge demand is negative when financial market and labor market is positively correlated, which means intuitively the investor shorts more or holds less of the stock overall, since income replaces holdings of the stock. As to the absolute amount of stock shorting, three factors have effects on in additional to ρ: the riskless rate $r$, the remaining time to the investment horizon $(T-t)$, and the ratio of volatility of labor market and financial market under the corresponding switching regime $\delta_i/\sigma_i$, i.e. the hedge demand is larger in magnitude if the risk rate is higher, volatility of income larger, volatility of stock smaller, and the investment horizon is further away.

Overall the optimal portfolio choice including the Merton component plus the hedge demand is determined by the signs of the risk premium and the correlation between financial market and labor market. If the risk premium and correlation have opposite signs, the optimal portfolio in the risky asset $\hat{\pi}$ is decreasing in magnitude with the investor's risk aversion, increasing in magnitude with the ratio of volatility of labor market and financial market. That is, if the correlation between stock returns and the income flow is negative, and risk premium positive, the optimal portfolio in the stock $\hat{\pi}$ is positive, decreasing with the investor's risk aversion, the volatility of stock market and increasing in magnitude with the volatility of labor market. Likewise, if the correlation between stock returns and the income flow is positive, and risk premium negative, the optimal portfolio in the stock $\hat{\pi}$ is negative, decreasing with the investor's risk aversion, the volatility of stock market and increasing in magnitude with the volatility of labor market.

On the contrary, if correlation and risk premium are of the same sign, the Merton component and the hedge demand will have opposite signs. In this case, the sign of the optimal portfolio could be either positive or negative and the investor may prefer to be either long or short the risky asset. For example, if both correlation and risk premium positive, the Merton component will be positive while the hedge term has negative sign. Thus overall the sign of optimal portfolio(longing or shorting the risky asset by the investor) is determined by the relative strength of the Merton component and the hedge demand. By analogy similar results can be reached if both correlation and risk premium are negative.

Furthermore, the regime switching plays an indispensable role in the portfolio choice. On the one hand, the switching parameters drive the risky asset's excess return and volatility, which have an impact on the Merton component of the optimal portfolio. On the other hand, as an important factor influencing the hedge part of the optimal portfolio, the ratio of volatility of labor market and financial market varies with different regimes.

## 5.Concluding remarks

In this paper, we set up a flexible and analytically tractable dynamic optimization model with dual Markov modulated diffusion processes, while the underlying Markov chains are not the same but distinct from each other. As an application and illustration of our model, we consider the optimal portfolio problem of an investor who wants to maximize his expected discounted utility of terminal wealth, which changes with the regime-switching financial market and labor market. After combining two single continuous time Markov chains into a composite Markov chain, we stated the corresponding HJB equations and obtain explicit solutions for the optimal control under some reasonable specifications and derived the corresponding value function up to an expectation over the Markov chain. We believe these establishment and solutions of stochastic control problems with dual Markov modulated diffusion processes constitute main contributions of this paper. We found out that the optimal portfolio consists of two parts: the first one corresponds to the solution of the Merton's classic mean-variance optimization problem and the second one corrects it for the additional risk coming from the stochastic labor income. We showed that both parts of the optimal portfolio depend on the state of the Markov chain. We observed that the investor hold a higher proportion of his wealth in the risky asset during a bull market regime than during a bear market regime, but the difference of the optimal stock position between the two states in absolute terms is reduced when the investor's risk tolerance goes down and enlarged with rising risk tolerance. The Markov chains have impact on the hedging part of the optimal policy in absolute terms via relative volatility of the stock market and the labor market under different regimes. However, the sign of the hedging term ,i.e. long or short of the risky asset out of consideration of hedging, is determined by the correlation between the Brownian motions for the stock price and the stochastic income processes.

**Appendix**

**Proof of Theorem 3.1**

*proof.* Let a and b satisfy $0 < a < |X(t)| < b < +\infty$. We define $\tau_a := \inf\{t \geq 0: |X(t)| = a\}$, $\tau_b := \inf\{t \geq 0: |X(t)| = b\}$, and $\tau := \tau_a \wedge \tau_b$. For $\forall t, \hat{t} \in [0, T]$, using the Itô formula for Markov-modulated processes(see, for instance, [29] [3] or [2]), we obtain

$$V(\hat{t} \wedge \tau, X(\hat{t} \wedge \tau), \xi(\hat{t} \wedge \tau)) - V(t, x, i)$$

$$= \int_t^{\hat{t} \wedge \tau} \left[ L_{\xi(s)} V(s, X(s), \xi(s)) \right.$$

$$\left. + \sum_{j=1}^{S} q_{\xi(s)j}(V(s, X(s), j) - V(s, X(s), \xi(s))) \right] ds$$

$$+ \int_t^{\hat{t} \wedge \tau} \sum_{k=1}^{p} \sum_{l=1}^{m} g_{kl}(s, X(s), \xi(s)) \frac{\partial V}{\partial x_k}(s, X(s), \xi(s)) dB_s^l$$

$$+ (M_{\hat{t} \vee \tau}^V - M_t^V) \quad (A.1)$$

where $\{M_t^V\}_{0 \leq t \leq T}$ is a real-valued, square integrable martingale, so we easily know $E[M_{\hat{t} \vee \tau}^V] = E[M_t^V]$. (A.2)

Also note that u is admissible and that $V'(s, X(s), \xi(s))$ is bounded for every $s \in [0, t \wedge \tau]$ (for each $i \in S$, $V(., i)$ is concave and $\max\{V'(\tau_a, a, \xi(\tau_a)), V'(\tau_b, b, \xi(\tau_b))\} < M$, where M is a finite constant). Then,

$E\left[\int_t^{\hat{t}\wedge\tau} \sum_{k=1}^{p} \sum_{l=1}^{m} g_{kl}(s, X(s), \xi(s)) \frac{\partial V}{\partial x_k}(s, X(s), \xi(s)) dB_s^l\right] = 0$ (A.3) because

$$E\left|\int_t^{\hat{t}\wedge\tau} \sum_{k=1}^{p} \sum_{l=1}^{m} g_{kl}(s, X(s), \xi(s)) \frac{\partial V}{\partial x_k}(s, X(s), \xi(s)) dB_s^l\right|^2 < +\infty$$

for all $t \in [0, T]$ by assumption 2.1.

Taking the expectation with respect to (A.1), together with (A.2) and (A.3), we have

$$E[V(\hat{t}, X(\hat{t}), \xi(\hat{t}))$$
$$= V(t, x, i)$$
$$+ E\{\int_t^{\hat{t}\wedge\tau}\left[L_{\xi(s)}V(s, X(s), \xi(s))\right.$$
$$\left.+ \sum_{j=1}^{S} q_{\xi(s)j}(V(s, X(s), j) - V(s, X(s), \xi(s)))\right]ds\}$$

(A.4)

Combining the principle of optimality(3.9) and (A.4), we get

$$0 = \sup_{u\in\mathcal{U}} E\{\int_t^{\hat{t}\wedge\tau} \Phi(s, X(s), u(s), \xi(s))ds$$
$$+ \int_t^{\hat{t}\wedge\tau}\left[L_{\xi(s)}V(s, X(s), \xi(s))\right.$$
$$\left.+ \sum_{j=1}^{S} q_{\xi(s)j}(V(s, X(s), j) - V(s, X(s), \xi(s)))\right]ds\}$$

Letting $a \downarrow 0$ and $b \uparrow +\infty$, we get $\tau \to T$. Therefore, applying the Monotone Convergence Theorem, we have that

$$0 = \sup_{u\in\mathcal{U}} E\{\int_t^{\hat{t}\wedge\tau} \Phi(s, X(s), u(s), \xi(s))ds$$
$$+ \int_t^{\hat{t}}\left[L_{\xi(s)}V(s, X(s), \xi(s))\right.$$
$$\left.+ \sum_{j=1}^{S} q_{\xi(s)j}(V(s, X(s), j) - V(s, X(s), \xi(s)))\right]ds\}$$

(A.5)

Letting $\hat{t} \to t$, we obtain the equation (3.11). The terminal condition of (3.11) holds obviously. Thus, the proof is done. ∎

**Proof of Proposition 4.3**

*Proof.* Firstly given each state $i \in \mathcal{S}$ and function $g([0,T] \times \mathbb{R}^N; \mathbb{R})$, define an operator $\mathcal{L}$ by

$$\mathcal{L}g(t,y,i) := g_y \mu(y,i,t) + \frac{1}{2} g_{yy} \delta(y,i,t)^2 + \sum_j q_{ij}(g(t,y,j) - g(t,y,i))$$

Meanwhile in order to simplify notation, let

$$c(s, X(s), \xi(s)) = -(e^{r(T-s)} \gamma Y_s + \frac{(\alpha_{\xi(s)} - r)^2}{2\sigma_{\xi(s)}^2})$$

$$Z(t) = exp\left(\int_0^t c(s, X(s), \xi(s)) ds\right)$$

According to the generalized Ito's lemma for Markov-modulated diffusion(see, e.g. [29][3][2]),for each state $\xi(t) \in \mathcal{S}$ and function $g([0,T] \times \mathbb{R}^N; \mathbb{R})$, we have

$$g(t, X(t), \xi(t)) - g(0, X(0), \xi(0)) = \int_0^t \left[(\frac{\partial}{\partial t} + \mathcal{L}) g(s, X(s), \xi(s))\right] ds + M_t^g$$

where $M_t^g$ is a local martingale.

and

$$E^{x,i} g(X(t), \xi(t)) - g(x,i) = E^{x,i} \int_0^t \left[(\frac{\partial}{\partial t} + \mathcal{L}) g(s, X(s), \xi(s))\right] ds$$

Applying Ito's formula to the switching process $g(Y(s), s, \xi(s)) Z(s)$, we have

$$E^{x,i} g(X(T), \xi(T)) Z(T) - g(x,i) = E^{x,i} \int_0^t \left[\left(\frac{\partial}{\partial t} + \mathcal{L}\right) g(s, X(s), \xi(s)) Z(s)\right] ds =$$

$$E^{x,i} \int_0^t Z(s) \left[g(s, X(s), \xi(s)) c(s, X(s), \xi(s)) + \mathcal{L}g(s, X(s), \xi(s)) + \frac{\partial g(s, X(s), \xi(s))}{\partial t}\right] ds = 0$$

Noting the boundary conditions, we obtain

$$g(t,y,i) = E^{x,i} g(X(T), \xi(T)) Z(T) = E^{x,i} Z(T) = E^{y,i}[exp\left(-\int_t^T (e^{r(T-s)} \gamma Y_s + \frac{(\alpha_{\xi(s)} - r)^2}{2\sigma_{\xi(s)}^2}) ds\right)]$$

substitute into (4.14) and (4.15),we get (4.21) and (4.22) respectively. The proof is completed. ∎